\documentclass[12pt]{article}

\usepackage{amsmath}
\usepackage{multirow} 

\usepackage[cp850]{inputenc}
\usepackage{latexsym,amsfonts}

\usepackage{graphicx}
\usepackage{subcaption}
\usepackage{pdflscape}

\usepackage{caption}

\usepackage{tikz-cd} 
\usepackage{tikz}

\usepackage{enumitem} 
\usepackage{authblk}

\newtheorem{Thm}{Theorem}[section]
\newtheorem{Note}{Note}

\author[]{Mustafa Gezek}
\affil[]{Department of Mathematics, Tekirdag Namik Kemal University, 
Tekirdag, Turkey 59030, \\ mgezek@nku.edu.tr}

\title{\bf New unitals in projective planes of order 16}
\title{\bf New unitals in projective planes of order 16}

\begin{document}
\maketitle

\begin{abstract}

In this study we performed a computer search for unitals in planes of order 16. Some new unitals were found and we show that some unitals can be embedded in two or more different planes. 
\medskip


\end{abstract}

\vspace{5mm}

\section{Introduction}

We assume familiarity with the basic facts
from combinatorial design theory and finite geometries 
%\cite{AK,BJL,CRC,Hir}.
\cite{BJL,CRC,Hir}.

A t-$(v,k,\lambda)$ design (or shortly, t-design) is a pair $D$=$\{ X, B \}$ satisfying the following conditions: 

\begin{enumerate}[label=(\roman*)]
\item $X$ has size $v$, called the {\it points} of $D$,
\item $B$ is a collection of subsets
of $X$ of size $k$, called the {\it blocks} of $D$,
\item Every $t$ points
 appear together in
exactly $\lambda$ blocks. % \cite{BJL}, \cite{CRC}.
\end{enumerate}

% (1) $X$ has size $v$, called the {\it points} of $D$; (2)$B$ is a collection of subsets of $X$ of size $k$, called the {\it blocks} of $D$; (3) every $t$ points appear together in exactly $\lambda$ blocks. % \cite{BJL}, \cite{CRC}.

Every point of a
$2$-$(v,k,\lambda)$ design is contained in
$$r=\lambda\frac{(v-1)}{(k-1)}$$ blocks, and the total number
of blocks is $$b=\frac{v}{k}r.$$ 
Here $r$ is called the {\it replication number} of the 2-design.

The incidence matrix of a $2$-$(v,k,\lambda)$ design $D$ is a $b \times v$ matrix with the property that the $(i,j)$ entry is 1 if $i$th block contains the $j$th point, 0 otherwise. For a prime  $p$, the rank
of the incidence matrix of a design $D$ over a finite field of characteristic $p$ is called the {\it $p$-rank} of $D$.

Two designs $D$ and $D'$ are called {\it isomorphic} if there is a bijection between their
point sets that maps every block of $D$ to a block of $D'$. An isomotphism of a design $D$ with itself is called an {\it automorphism} of $D$.

The {\it dual} design of a design $D$ is defined to be the design having the incidence matrix as the transpose of the incidence matrix of $D$, denoted by $D^\perp$. In other words, the dual design of a design $D$ has point set as the block set of $D$, and the block set as the point set of $D$.

If the number of points of a design $D$ is equal to the number of blocks, then $D$ is called a {\it symmetric} design.

A $2$-$(v,k,1)$ design is called a {\it Steiner} design. 
A projective plane of order $q$ is a symmetric
Steiner 2-$(q^2 + q +1, q+1, 1)$ design with $q\ge 2$.Desarguesian plane $PG(2,p^t)$
of order $q=p^t$, where $p$ is prime and $t\ge 1$, has as points
the 1-dimensional subspaces of the 3-dimensional vector space $V_3$
over the finite field of order $q$, and as blocks (or {\it lines}), the
2-dimensional subspaces of $V_3$.

Point sets in projective planes with two line intersections
have been studied in finite geometry a lot.
{\it Maximal arcs} and {\it unitals} are examples of such sets.

A {\it unital} in a projective plane $\Pi$ of order $q^2$ is a set $U$ of $q^3+1$ points that meets every line of the plane in either one or $q+1$ points. If $U$ is a unital in $\Pi$, the lines that meet $U$ in one point form a unital $U^\perp$ in $\Pi^\perp$. The point set of a unital with the line intersections of size $q+1$ form a $2-(q^3+1,q+1,1)$ design, called {\it unital design associated $U$.}

In this study, the results of a computer search for unitals in planes of order 16 are reported. Our computations confirm most of the data given in \cite{S1}, but provides more sets of unitals in some specific group orders.

In Section \ref{sec2}, we provide a brief history on unitals and unitals in planes of order 16 are discussed in details. The number of known unitals in these planes with their automorphism group order are given in Table \ref{tab1}. Table \ref{tab1} shows that at least 20 new unitals were discovered by the algorithm we have.  

In Section \ref{sec3}, details of $2-(65,5,1)$ designs associated to unitals given in Section \ref{sec2} are provided: the orders of the automorphism group, the 5-rank and the number of parallel classes of designs associated to unitals and dual unitals are computed. The bound on the number of non-isomorphic unital designs embedded in planes of order 16 is improved.

In the last section, we list explicitly the point sets of the unitals given in Table \ref{tab2}. 28 of these 93 unitals given in Section \ref{sec4} already listed in \cite{ST1}.

\section{Unitals in planes of order 16}
\label{sec2}
13 projiective planes of order 16 are known to exist, four of which are self-dual, and the rest have non-isomorphic duals. The names of the planes are in accordance with \cite{PRS}: BBH1, BBH2, BBS4, DEMP, DSFP, HALL, JOHN, JOWK, LMRH, MATH, SEMI2, SEMI4 and PG(2,16). 

Stoichev and Tonchev \cite{ST1} performed a non-exhaustive search for unitals in the known planes of order 16 and they found 38 such examples, which are explicitly listed in their paper. Later, Stoichev \cite{S1} claimed that the number of known unitals in planes of order 16 is 116, and he gave the order of their automorphism group in his paper, but the point sets of the unitals neither were given in the paper nor provided online at the time of writing this paper. 

In 2011, Kr\v{c}adinac and Smoljak \cite{KS} claimed that they found three more unitals, one in LMRH plane with group order 16 and two in HALL plane with group orders 16 and 80, respectively. Most probably these unitals were already known by Stoichev by that time since \cite{S1} claims that these planes already have the unitals with these group orders. When we checked isomorphisms between these unitals with the ones we found, we see that these are listed in Table \ref{tab2}: the unital $\#1$ in LMRH plane, and the unitals $\#1$ and $\#4$ in HALL planes, respectively.

The specific line sets of the known projective planes of order 16
that we are using in this paper were graciously provided
by Gordon F. Royle, and are available online at 
\begin{verbatim}
http://pages.mtu.edu/~tonchev/planesOForder16.txt
\end{verbatim}

We performed a non-exhaustive search for unitals in the planes of order 16 and we found some new unitals. The algorithm we used is based on unions of orbits of appropriate subgroups of the
the automorphism group of the associated plane. Only subgroups of order 12, 16 and 20 were used in this search.

Table \ref{tab1} gives the number of unitals we found and the ones already known in planes of order 16. Column 1 lists the names of the planes, Columns 2-16 gives the number of known unitals with the specific automorphism group orders, and the last column shows the total number of known unitals in each plane.

\begin{table}[h!]
\caption{Number of unitals in projective planes of order 16}
\medskip
\centering\renewcommand{\arraystretch}{1.2}
\scalebox{.6}{
\begin{tabular}{|c|c|c|c|c|c|c|c|c|c|c|c|c|c|c|c|c|c|c|}
\hline
 Plane/$|Aut(Unital)|$ & 4 & 8 & 12 & 16  & 20 & 24 & 32 & 48 & 64 & 80 & 100 & 128 & 192 & 1200 & 768 & 249600 & TOTAL \\ \hline
BBH1   &   & 1 &   & 4 &   &   & 2 &   &   &   &  &  &  &  &  &  & 7\\  \hline
BBH2   & 1 & 1 &   & 8 & 1 &   & 2 &   &   & 1 &  &  &  &  &  &  & 14 \\  \hline
BBS4   &   & 1 & 2 &   &   & 1 &   &   &   &   &  &  &  &  &  &  &4 \\  \hline
DEMP   &   &   & 1 & 1 &   & 1 &   & 1 &   &   &  &  &  &  &   &  &4\\  \hline
DSFP   &   &   & 1 &   &   & 1 &   &   &   &   &  &  &  &  &  &  & 2 \\  \hline
HALL   &   &   &   & 1 &   &   & 1 & 1 &   & 1 & 1&  &  & 1&  &  & 6\\  \hline 
LMRH   &   &  &   & 1 &   &   & 1 &   &   &   &  &  &  &  &  &  & 2\\  \hline
JOHN   &   & 1 &   & 4 &   & 1 & 2 & 1 &   &   &  &  &  &  &  &  & 9\\  \hline
JOWK   & 1 & 1 & 1 & 1 &   & 1 & 1 & 1 &   &   &  &  &  &  &  &  & 7\\  \hline
MATH   &   &   &   & 6 &   &   & 1 &   & 3 &   &  & 6&  &  &   & &16\\  \hline
SEMI2  &   &   &  & 1 &   &   & 2 & 3 & 12&   &  &  & 3&  &  &  &21\\  \hline
SEMI4  &   & 1  & 2 &   &   &   &   & 1 & 1 &   &  & 3& 1&  &  &  &9 \\  \hline
PG(2,16)&   &   &   &   &   &   &   &   &  &  & &  & &  & 1 & 1 &  2 \\  \hline
% TOTAL.  &   &   &   &   &   &   &   &   &  &  &  &  & \\  \hline

\end{tabular}}
\label{tab1}
\end{table}

When we compare the data given in Table \ref{tab1} with the data given in \cite{S1}, we see that there are {\it at least} 20 new unitals in planes of order 16: {\it one} in BBH1, {\it three} in BBS4, {\it one} in DEMP, {\it two} in JOHN, {\it two} in MATH and {\it twelve} in SEMI2.  

Previously, it was known that some unital 2-(28,4,1) can be embedded in a plane of order 9 and its dual,
and that plane was not self-dual. Our computations show that the unital $\#6$ in SEMI4 is a unital not only in SEMI4, but also in JOWK and HALL planes, and the unital $\#5$ in SEMI4 is also a unital in HALL plane. This observation shows that some unitals can be embedded in two or more different planes.

%The point sets of all found unitals with pur algorithm can be found in Section \ref{sec4}.

\section{Unital Designs}
\label{sec3}

Table \ref{tab2} contains some data related to unital designs found with our algorithm. Column 3 lists the order of the automorphism group of each unital, Column 4 gives the 5-rank of each designs, and Column 5 provides the number of parallel classes of each design of unital and design of the dual unital.
The last column indicates if the unital design is isomorphic to any other designs, including designs of dual unitals. 

%The order of the automorphism groups, 5-ranks and number of parallel classes of each unital designs were computed, and they are listed in Column 3 to Column 5, respectively. 

In Table \ref{tab2} a $^*$ indicates that the unital was already known before and listed in \cite{ST1}, implying that 28 unitals listed in Table \ref{tab2} were already known before. \cite{ST1} also contains eight unitals having group order $<9$ and two unitals in PG(2,16), the details of these unital designs are listed in Table \ref{tab3}, in which the unital numbers are as the continuation of Table \ref{tab2}.

The total number of designs given in Table \ref{tab2} and Table \ref{tab3} is 103, and they are pairwise non-isomorphic. Considering the designs of dual unitals, we see that the number of 
unital designs embedded in planes of order 16 is at least 168, this number previously was 73 \cite{KNP}.

\begin{Thm}
The number of non-isomorphic unital designs embedded in projective planes of order 16 is $\geq$ 168. 
\end{Thm}{}

\begin{table}[h!]
\caption{Designs associated to unitals in planes of order 16}
\medskip
\centering\renewcommand{\arraystretch}{1.2}
\scalebox{.67}{
\begin{tabular}{|c|c|c|c|c|c|}
\hline
 Plane & Unital $\#$  & $|Aut(Unital)|$  & $5-$rank  & $\#$ of par. clas. & Isomorphic to?  \\ \hline
BBH1 & 1 & 16  & 63  & 34/34 & BBH$1^\perp.3$ \\ 
     & 2 & 16  & 62  & 90/138 & BBH$1^\perp.4$ \\ 
     & 3 & 16  & 63  & 34/34 & BBH$1^\perp.1$\\ 
     & 4$^*$ & 16  & 62  & 138/90 & BBH$1^\perp.2$\\ 
     & 5$^*$ & 32  & 64  & 32/32 & BBH$1^\perp.5$\\ 
     & 6 & 32  & 64  & 24/24 & BBH$1^\perp.6$\\      \hline
BBH2 & 1 & 16  & 64  & 28/32 &  \\ 
     & 2 & 16  & 64  & 36/24 &   \\ 
     & 3 & 16  & 64  & 32/36 &  \\ 
     & 4 & 16  & 64  & 36/36 &  \\ 
     & 5 & 16  & 64  & 28/32 &  \\ 
     & 6 & 16  & 63  & 38/34 &  \\ 
     & 7 & 16  & 63  & 50/62 &  \\ 
     & 8$^*$ & 16  & 63  & 34/38 &  \\ 
     & 9$^*$ & 20  & 60  & 119/64 &  \\ 
     & 10$^*$ & 32  & 64  & 24/16 &  \\ 
     & 11 & 32  & 64  & 48/40 &  \\ 
     & 12$^*$ & 80  & 59  & 174/154 &  \\      \hline
BBS4 & 1 & 12  & 64  & 4/4 &  \\ 
     & 2 & 12  & 63  & 34/34 &   \\ 
     & 3 & 24  & 63  & 38/34 &  \\ \hline
DEMP & 1$^*$ & 12  & 64  & 4/4 &  \\ 
     & 2 & 16  & 64  & 10/14 &  \\ 
     & 3$^*$ & 24  & 64  & 10/10 &  \\ 
     & 4 & 48  & 64  & 10/10 &  \\ \hline
DSFP & 1$^*$ & 12  & 63  & 34/34 &  \\ 
     & 2$^*$ & 24  & 63  & 38/34 &   \\ \hline
HALL & 1 & 16  & 64  & 20/16 &  \\ 
     & 2$^*$ & 32  & 64  & 48/48 &   \\ 
     & 3$^*$ & 48  & 63  & 34/34 &  \\ 
     & 4 & 80  & 60  & 224/184 &  \\ 
     & 5$^*$ & 100  & 58  & 264/114 &  \\ 
     & 6$^*$ & 1200  & 54  & 1094/1094 & HALL$^\perp.6$ \\ \hline
JOHN & 1 & 16  & 64  & 28/24 &  \\ 
     & 2 & 16  & 63  & 34/42 &   \\ 
     & 3$^*$ & 16  & 63  & 42/42 &  \\ 
     & 4 & 16  & 63  & 46/42 &  \\ 
     & 5$^*$ & 24  & 63  & 34/34 &  \\ 
     & 6$^*$ & 32  & 64  & 32/32 &  \\ 
     & 7 & 32  & 64  & 16/16 &  \\ 
     & 8$^*$ & 48  & 63  & 34/34 &  \\      \hline
JOWK & 1 & 12  & 64  & 4/4 &  \\ 
     & 2 & 16  & 63  & 38/46 &   \\ 
     & 3$^*$ & 24  & 63  & 34/34 &  \\ 
     & 4 & 32  & 64  & 16/16 &  \\ 
     & 5$^*$ & 48  & 63  & 34/34 &  \\  \hline

\end{tabular}}
\label{tab2}
\end{table}

\begin{table}[]
\ContinuedFloat  %% <-- new
\caption{Designs associated to unitals in planes of order 16 (continued)}
\medskip
\centering\renewcommand{\arraystretch}{1.2}
\scalebox{.65}{
\begin{tabular}{|c|c|c|c|c|c|}
\hline
 Plane & Unital $\#$  & $|Aut(Unital)|$  & $5-$rank  & $\#$ of par. clas. & Isomorphic to?  \\ \hline
LMRH & 1 & 16  & 64  & 16/16 &  \\ 
     & 2$^*$ & 32  & 64  & 16/16 &   \\ \hline
MATH & 1 & 16  & 64  & 32/28 &  \\ 
     & 2$^*$ & 16  & 64  & 4/4 &   \\ 
     & 3 & 16  & 64  & 0/0 &  \\ 
     & 4 & 16  & 64  & 0/0 &  \\ 
     & 5 & 16  & 64  & 0/0 &  \\ 
     & 6 & 16  & 64  & 22/22 &  \\ 
     & 7$^*$ & 32  & 64  & 16/16 &  \\ 
     & 8$^*$ & 64  & 64  & 16/16 &  \\ 
     & 9 & 64  & 64  & 16/16 &  \\ 
     & 10 & 64  & 64  & 16/16 &  \\ 
     & 11 & 128  & 64  & 16/16 &  \\ 
     & 12$^*$ & 128  & 64  & 16/16 &  \\ 
     & 13 & 128  & 64  & 80/80 &  \\ 
     & 14 & 128  & 64  & 16/16 &  \\ 
     & 15 & 128  & 64  & 16/16 &  \\ 
     & 16 & 128  & 64  & 16/16 &  \\      \hline
SEMI2 & 1 & 16  & 64  & 4/4 & SEMI$2^\perp.1$ \\ 
     & 2 & 32  & 64  & 16/16 & SEMI$2^\perp.3$  \\ 
     & 3 & 32  & 64  & 16/16 & SEMI$2^\perp.2$ \\ 
     & 4 & 48  & 64  & 4/4 & SEMI$2^\perp.5$ \\ 
     & 5 & 48  & 64  & 4/4 & SEMI$2^\perp.4$ \\ 
     & 6$^*$ & 48  & 64  & 4/4 & SEMI$2^\perp.6$ \\ 
     & 7 & 64  & 64  & 16/16 & SEMI$2^\perp.7$ \\ 
     & 8 & 64  & 64  & 32/32 & SEMI$2^\perp.9$ \\ 
     & 9$^*$ & 64  & 64  & 32/32 & SEMI$2^\perp.8$ \\ 
     & 10 & 64  & 64  & 16/16 & SEMI$2^\perp.10$ \\ 
     & 11 & 64  & 64  & 16/16 & SEMI$2^\perp.12$ \\ 
     & 12 & 64  & 64  & 16/16 & SEMI$2^\perp.11$ \\ 
     & 13 & 64  & 64  & 16/16 & SEMI$2^\perp.17$ \\ 
     & 14 & 64  & 64  & 16/16 & SEMI$2^\perp.18$ \\ 
     & 15 & 64  & 64  & 16/16 & SEMI$2^\perp.15$ \\ 
     & 16 & 64  & 64  & 16/16 & SEMI$2^\perp.16$ \\ 
     & 17 & 64  & 64  & 16/16 & SEMI$2^\perp.13$ \\ 
     & 18 & 64  & 64  & 16/16 & SEMI$2^\perp.14$ \\ 
     & 19 & 192  & 64  & 16/16 & SEMI$2^\perp.20$ \\ 
     & 20 & 192  & 64  & 16/16 & SEMI$2^\perp.19$ \\ 
     & 21$^*$ & 192  & 64  & 16/16 & SEMI$2^\perp.21$ \\      \hline     
SEMI4 & 1 & 12  & 64  & 16/16 & SEMI$4^\perp.2$ \\ 
     & 2 & 12  & 64  & 16/16 & SEMI$4^\perp.1$  \\ 
     & 3 & 48  & 64  & 4/4 & SEMI$4^\perp.3$ \\ 
     & 4 & 64  & 64  & 16/16 & SEMI$4^\perp.4$ \\ 
     & 5 & 128  & 64  & 16/16 & SEMI$4^\perp.7$ \\ 
     & 6 & 128  & 64  & 16/16 & SEMI$4^\perp.6$ \\ 
     & 7 & 128  & 64  & 16/16 & SEMI$4^\perp.5$ \\ 
     & 8 & 192  & 64  & 16/16 & SEMI$4^\perp.8$ \\      \hline
%PG(2,16) & 1 & 768  &   &  &  \\ 
%     & 2 & 249600  &   &  &   \\  \hline
\end{tabular}}
\label{tab2}
\end{table}

\begin{table}[h!]
\caption{Designs associated to unitals in planes of order 16}
\medskip
\centering\renewcommand{\arraystretch}{1.2}
\scalebox{.67}{
\begin{tabular}{|c|c|c|c|c|c|}
\hline
 Plane & Unital $\#$  & $|Aut(Unital)|$  & $5-$rank  & $\#$ of par. clas. & Isomorphic to?  \\ \hline
BBH1 & 7 & 8  & 63  & 34/34 &  \\ \hline
BBH2 & 13 & 4  & 63  & 18/18 &  \\ 
     & 14 & 8  & 61  & 82/86 &  \\ \hline
BBS4 & 4 & 8  & 63  & 38/42 &  \\ \hline
JOHN & 9 & 8  & 63  & 34/36 &  \\ \hline
JOWK & 6 & 4  & 64  & 25/16 &  \\ 
     & 7 & 8  & 63  & 34/34 &  \\ \hline
SEMI4 & 9 & 8  & 64  & 20/24 &  \\ \hline
PG(2,16) & 1 & 768  & 64  & 16/16 & PG(2,16)$^\perp.1$  \\ 
 & 2 & 249600  & 52  & 4304/4304 & PG(2,16)$^\perp.2$ \\ \hline
\end{tabular}}
\label{tab3}
\end{table}

\newpage

\section{The point sets of unitals}
\label{sec4}
In this section we list only the point sets of unitals given in Table \ref{tab2}, the point sets of the rest can be found in \cite{ST1}.

\vspace{.03in}

{\bf BBH1:}

\vspace{.03in}

    $\{ 257, 262, 263, 265, 264, 82, 87, 96, 90, 18, 185, 20, 30, 184, 183, 186, 
    26, 33, 131,\\38, 99, 235, 144, 154, 55, 104, 233, 194, 242, 160, 49, 200, 
    244, 1, 4, 15, 68, 122, 2, 77, \\ 113, 76, 120, 224, 212, 67, 128, 215, 219, 
    37, 129, 43, 106, 232, 139, 147, 60, 102, 227, \\ 204, 245, 153, 58, 207, 253 
    \}$,

\vspace{.03in}
    
    $\{ 257, 262, 263, 265, 264, 82, 87, 96, 90, 18, 185, 20, 30, 184, 183, 186, 
    26, 40, 133, \\ 
    44, 109, 226, 143, 159, 62, 112, 239, 202, 249, 152, 61, 208, 
    246, 8, 14, 10, 74, 126, 11, \\ 72, 124, 73, 117, 217, 218, 80, 114, 210, 209, 
    42, 142, 45, 101, 229, 135, 158, 53, 107, \\ 231, 206, 251, 150, 56, 203, 243 
    \}$,

\vspace{.03in}
    
    $\{ 267, 21, 163, 187, 85, 8, 123, 77, 209, 48, 256, 196, 227, 104, 145, 135,
    52, 37, 251, \\198, 231, 102, 151, 136, 152, 105, 228, 206, 49, 59, 250, 41, 
    139, 40, 243, 203, 236, 106, \\153, 143, 156, 98, 230, 200, 63, 57, 244, 46, 
    140, 19, 161, 177, 180, 175, 88, 92, 30, 24, \\96, 90, 186, 166, 184, 168, 26 
    \}$,

\vspace{.03in}
    
    $\{ 267, 3, 122, 70, 213, 28, 171, 192, 84, 39, 245, 208, 237, 107, 155, 144,
    62, 37, 251, \\198, 231, 102, 151, 136, 152, 105, 228, 206, 49, 59, 250, 41, 
    139, 36, 255, 193, 229, 103, \\146, 134, 159, 112, 235, 199, 53, 55, 247, 43, 
    129, 19, 161, 177, 180, 175, 88, 92, 30, 24, \\96, 90, 186, 166, 184, 168, 26 
    \}$,

\vspace{.03in}
    
    $\{ 257, 83, 188, 92, 94, 179, 182, 189, 88, 5, 15, 215, 13, 213, 11, 211, 
    209, 66, 124, 75, \\72, 122, 123, 68, 125, 18, 162, 20, 30, 173, 165, 26, 167,
    49, 134, 149, 59, 64, 248, 143, \\140, 245, 150, 160, 56, 256, 129, 246, 148, 
    33, 102, 195, 39, 48, 240, 109, 105, 237, 206, \\194, 42, 227, 97, 226, 198 
    \}$,

\vspace{.03in}
    
    $\{ 257, 114, 222, 120, 119, 118, 224, 220, 218, 17, 192, 25, 23, 32, 178, 
    181, 187, 2, 73, \\3, 8, 6, 79, 70, 77, 82, 164, 170, 90, 87, 96, 166, 161, 
    51, 129, 159, 62, 52, 58, 246, 140, 143, \\248, 156, 153, 134, 146, 256, 245, 
    37, 197, 232, 44, 41, 47, 109, 201, 207, 105, 239, 225, \\202, 234, 102, 97 \}$.

\vspace{.15in}
    
{\bf BBH2:}

\vspace{.05in}

     $\{ 83, 53, 55, 141, 132, 242, 251, 145, 155, 34, 159, 231, 104, 193, 254, 52, 139, 18, 106, \\ 169, 89, 208, 41, 190, 234, 50,
    60, 143, 245, 152, 140, 256, 157, 97, 247, 197, 105, 109, 160, \\ 243, 255, 149, 201, 207, 102, 150, 244, 148, 202, 5, 9, 119,
    71, 215, 117, 113, 73, 7, 213, 80, \\114, 216, 214, 77, 15 \}$,

\vspace{.02in}
       
    $\{ 83, 53, 158, 252, 132, 242, 131, 155, 63, 1, 14, 123, 218, 125, 224, 75, 66, 217, 257, 221, \\223, 269, 270, 219, 259, 97,
    149, 102, 105, 109, 160, 148, 150, 26, 37, 91, 176, 86, 98, 239, \\103, 44, 191, 30, 177, 196, 205, 232, 163, 19, 198, 162, 
    185, 173, 101, 42, 110, 194, 81, 31, \\84, 236, 235, 48, 183 \}$,
 
\vspace{.02in}
       
    $\{ 83, 28, 29, 188, 90, 182, 82, 170, 164, 217, 257, 221, 223, 269, 270, 219, 259, 50, 248, \\151, 245, 152, 140, 61, 129, 
    97, 149, 102, 105, 109, 160, 148, 150, 18, 225, 95, 169, 89, 208, \\47, 204, 234, 180, 20, 190, 106, 112, 41, 171, 4, 71, 77,
    115, 79, 120, 5, 114, 78, 216, 7, 210, \\209, 215, 13, 117 \}$,
 
\vspace{.02in}
      
    $\{ 83, 36, 144, 43, 46, 40, 130, 135, 137, 21, 96, 179, 87, 189, 161, 22, 166, 54, 157, 246, \\147, 143, 60, 134, 256, 121, 
    260, 128, 126, 272, 267, 127, 262, 35, 51, 203, 107, 199, 138, \\153, 133, 62, 228, 38, 233, 253, 241, 146, 111, 19, 198, 
    162, 185, 173, 101, 42, 110, 194, 81, \\31, 84, 236, 235, 48, 183 \}$,

\vspace{.02in}
       
    $\{ 83, 55, 250, 141, 251, 154, 57, 145, 142, 217, 257, 221, 223, 269, 270, 219, 259, 97, 149, \\102, 105, 109, 160, 148, 150,
    12, 68, 124, 220, 122, 222, 72, 16, 34, 156, 200, 104, 193, 254, \\58, 249, 159, 238, 45, 231, 139, 136, 52, 100, 9, 11, 70, 
    113, 73, 213, 76, 211, 15, 212, 10, 214, \\119, 118, 80, 116 \}$,

\vspace{.02in}
       
    $\{ 32, 17, 23, 83, 25, 7, 73, 212, 115, 35, 43, 231, 111, 200, 202, 97, 240, 48, 95, 58, 108, \\206, 242, 89, 86, 254, 146, 
    235, 91, 53, 132, 133, 155, 37, 85, 255, 98, 208, 157, 81, 84, 59, \\135, 225, 93, 143, 60, 160, 256, 1, 67, 271, 66, 218, 
    76, 6, 127, 273, 272, 125, 217, 13, 114, \\270, 213 \}$,

\vspace{.02in}
       
    $\{ 32, 22, 87, 170, 182, 17, 23, 83, 25, 35, 43, 231, 200, 111, 202, 240, 97, 3, 263, 71, 127, \\221, 224, 270, 260, 114, 5, 
    65, 268, 125, 72, 215, 16, 31, 243, 152, 183, 81, 106, 149, 144, 61, \\143, 162, 59, 103, 112, 246, 98, 18, 62, 142, 163, 
    180, 42, 249, 133, 141, 131, 86, 159, 110, \\235, 132, 198 \}$,

\vspace{.02in}
        
    $\{ 32, 22, 87, 170, 182, 17, 23, 83, 25, 35, 43, 231, 200, 111, 202, 240, 97, 26, 156, 57, 171, \\191, 195, 138, 254, 55, 63,
    89, 51, 229, 108, 53, 33, 8, 259, 119, 126, 217, 12, 266, 265, 80, \\213, 74, 272, 68, 123, 15, 218, 19, 130, 245, 185, 84, 
    239, 56, 247, 129, 60, 173, 160, 234, \\232, 147, 225 \}$,
 
\vspace{.02in}
      
    $\{ 187, 172, 174, 192, 178, 22, 29, 188, 164, 166, 10, 126, 272, 66, 211, 
    26, 143, 226, 44, \\236, 35, 84, 142, 56, 62, 34, 208, 173, 165, 242, 160, 
    230, 132, 203, 101, 109, 159, 59, 89, \\246, 37, 227, 91, 134, 51, 27, 158, 
    199, 110, 106, 207, 31, 190, 50, 138, 243, 39, 177, 252, \\228, 105, 196, 148,
    153, 245 \}$,
 
\vspace{.02in}
      
    $\{ 83, 21, 164, 28, 87, 189, 161, 182, 90, 34, 159, 231, 104, 193, 254, 52, 139, 12, 75, 123, \\220, 122, 224, 14, 72, 18, 
    106, 169, 89, 208, 41, 190, 234, 2, 260, 127, 6, 3, 261, 267, 272, \\268, 126, 128, 8, 273, 262, 265, 121, 5, 9, 119, 71, 
    215, 117, 113, 73, 7, 213, 80, 114, 216, \\214, 77, 15 \}$,
 
\vspace{.02in}
    
    $\{ 83, 34, 159, 231, 104, 193, 254, 52, 139, 18, 106, 169, 89, 208, 41, 190, 234, 22, 82, 29, \\96, 179, 170, 188, 166, 1, 
    16, 124, 218, 125, 222, 66, 68, 97, 247, 197, 105, 109, 160, 243, \\255, 149, 201, 207, 102, 150, 244, 148, 202, 19, 195, 24,
    185, 173, 101, 33, 108, 236, 184,\\ 167, 84, 194, 229, 48, 93 \}$,
 
\vspace{.02in}
     
    $\{ 32, 17, 23, 83, 25, 7, 73, 212, 115, 35, 43, 231, 111, 200, 202, 97, 240, 8, 259, 119, 65, \\223, 12, 260, 261, 215, 78, 
    128, 258, 72, 220, 11, 122, 48, 95, 58, 108, 206, 242, 89, 86, 254, \\146, 235, 91, 53, 132, 133, 155, 19, 130, 245, 24, 27, 
    239, 148, 243, 50, 140, 31, 49, 196, 106,\\ 152, 47 \}$.

\vspace{.1in}
    
{\bf BBS4:}

\vspace{.03in}

  $\{ 53, 63, 259, 50, 60, 8, 58, 237, 47, 211, 17, 195, 253, 198, 214, 252, 
    236, 18, 21, 80, 57, \\88, 62, 68, 97, 103, 74, 127, 107, 66, 101, 225, 105, 
    190, 78, 232, 137, 202, 162, 149, 207, 76, \\125, 205, 99, 159, 86, 227, 186, 
    230, 145, 184, 204, 2, 5, 181, 41, 173, 46, 169, 140, 158, 163, \\134, 130 \}$,

\vspace{.02in}
         
    $\{ 53, 63, 50, 60, 259, 72, 146, 122, 81, 194, 176, 135, 111, 256, 215, 185,
    233, 5, 129, 59, \\20, 187, 191, 152, 46, 133, 174, 170, 148, 69, 230, 123, 
    84, 183, 220, 243, 110, 137, 162, 205,\\ 160, 8, 245, 58, 17, 206, 203, 228, 
    47, 244, 219, 222, 229, 10, 89, 56, 31, 71, 103, 80, 33, 121,\\ 82, 114, 112 
    \}$,
     
\vspace{.02in}
     
        $\{ 263, 273, 264, 271, 50, 4, 35, 62, 21, 32, 29, 54, 43, 12, 34, 9, 55, 71,
    177, 121, 82, 201, \\143, 168, 112, 154, 247, 224, 226, 5, 37, 59, 20, 31, 27,
    52, 46, 14, 33, 10, 56, 68, 187, 126, 85, \\200, 133, 174, 107, 148, 250, 209,
    239, 66, 183, 128, 87, 202, 137, 162, 105, 160, 248, 223, \\225 \}$.

\vspace{.1in}
    
{\bf DEMP:}

\vspace{.03in}

      $\{ 257, 46, 110, 94, 190, 8, 163, 124, 24, 58, 131, 69, 196, 236, 212, 154, 
    245, 2, 162, 127, \\18, 63, 130, 70, 198, 239, 214, 159, 246, 13, 173, 121, 
    29, 57, 141, 75, 203, 233, 219, 153, 251, \\ 1, 161, 128, 17, 64, 129, 71, 199,
    240, 215, 160, 247, 10, 164, 117, 26, 56, 132, 76, 195, 229, \\211, 152, 252 
    \}$, 
    
\vspace{.02in}

    $\{ 273, 21, 24, 31, 32, 1, 2, 10, 4, 51, 158, 54, 55, 60, 153, 155, 157, 75,
    131, 237, 123, 78,\\ 214, 167, 140, 253, 233, 195, 126, 204, 215, 166, 249, 
    72, 143, 240, 128, 80, 213, 165, 133,\\ 248, 232, 197, 120, 207, 223, 175, 
    256, 35, 46, 86, 183, 44, 89, 190, 43, 108, 87, 109, 182, 105, \\93, 187, 99 
    \}$,
    
\vspace{.02in}
     
    $\{ 80, 224, 257, 16, 192, 23, 225, 197, 106, 94, 175, 157, 38, 249, 136, 
    115, 50, 25, 237, 202, \\101, 88, 166, 145, 47, 247, 142, 114, 51, 20, 236, 
    196, 100, 91, 172, 156, 44, 244, 139, 123, 59, \\17, 238, 207, 109, 87, 168, 
    147, 34, 246, 133, 122, 57, 22, 227, 205, 111, 85, 162, 158, 40, 241, \\135, 
    121, 58 \}$,
    
\vspace{.02in}
     
        $\{ 257, 11, 59, 27, 155, 46, 110, 94, 190, 77, 217, 237, 125, 253, 137, 201,
    169, 67, 220, 236, \\119, 243, 131, 198, 204, 246, 124, 167, 231, 163, 70, 
    214, 135, 66, 223, 234, 113, 242, 136, 197,\\ 207, 244, 122, 176, 225, 168, 
    68, 213, 144, 33, 112, 84, 178, 81, 181, 88, 192, 106, 36, 47, 98, \\101, 186,
    40, 95 \}$.

\vspace{.1in}
    
{\bf DSFP:}

\vspace{.03in}

    $\{ 266, 273, 257, 260, 258, 41, 224, 187, 94, 119, 207, 136, 109, 70, 252, 
    165, 227, 25, 64, \\27, 30, 151, 63, 56, 29, 150, 156, 53, 147, 44, 221, 179, 
    87, 117, 206, 139, 102, 72, 255, 169, \\240, 34, 212, 178, 82, 122, 196, 132, 
    98, 74, 250, 164, 234, 47, 214, 192, 85, 121, 199, 131, \\104, 75, 254, 172, 
    237 \}$,
    
\vspace{.02in}
 
         $\{ 1, 266, 26, 50, 148, 36, 38, 180, 84, 150, 183, 83, 108, 100, 151, 147, 
    156, 48, 198, 191,\\ 88, 123, 215, 163, 140, 101, 73, 253, 238, 53, 242, 56, 
    63, 79, 226, 114, 66, 64, 128, 229, 248, \\21, 234, 24, 31, 72, 250, 74, 122, 
    32, 117, 240, 255, 3, 81, 12, 6, 103, 97, 33, 177, 7, 86, 188, \\35 \}$.

\vspace{.1in}
    
{\bf HALL:}

\vspace{.03in}

   $\{ 270, 70, 156, 199, 83, 42, 49, 210, 244, 43, 72, 249, 62, 221, 208, 149, 
    95, 3, 87, 119, 108, \\231, 166, 131, 150, 195, 6, 19, 134, 236, 188, 103, 76,
    13, 64, 222, 105, 235, 251, 16, 40, 255, \\133, 41, 142, 104, 61, 239, 213, 9,
    159, 203, 109, 238, 78, 63, 85, 80, 216, 93, 139, 37, 153, \\256, 200 \}$,
   
\vspace{.02in}
 
        $\{ 260, 100, 111, 116, 157, 220, 156, 221, 127, 56, 183, 146, 238, 105, 202,
    115, 213, 49, \\182, 187, 230, 64, 193, 208, 235, 6, 241, 79, 96, 29, 171, 36,
    140, 5, 245, 145, 83, 249, 112, \\163, 9, 150, 107, 118, 162, 219, 82, 113, 
    224, 12, 253, 229, 95, 244, 227, 164, 15, 197, 195, \\57, 173, 50, 92, 185, 
    178 \}$,
   
\vspace{.02in}
 
       $\{ 266, 42, 81, 178, 100, 257, 260, 258, 273, 15, 243, 160, 21, 119, 76, 
    230, 56, 3, 35, 196, \\194, 193, 61, 29, 13, 150, 60, 157, 23, 86, 188, 103, 
    202, 10, 211, 78, 126, 254, 96, 47, 101, \\145, 50, 184, 20, 214, 220, 215, 
    238, 41, 225, 72, 117, 255, 176, 175, 165, 89, 185, 168, 105, \\250, 116, 66, 
    240 \}$,
    
\vspace{.02in}
 
       $\{ 270, 70, 156, 199, 83, 42, 49, 210, 244, 43, 72, 249, 62, 221, 208, 149, 
    95, 4, 212, 82, 98,\\ 225, 148, 66, 242, 33, 81, 202, 138, 196, 65, 154, 58, 
    13, 64, 222, 105, 235, 251, 16, 40, 255, \\133, 41, 142, 104, 61, 239, 213, 8,
    125, 152, 112, 229, 96, 110, 169, 27, 237, 69, 143, 11, 207,\\ 137, 190 \}$,
   
\vspace{.02in}
 
        $\{ 46, 245, 231, 195, 141, 261, 263, 272, 262, 265, 12, 65, 22, 122, 104, 3,
    225, 21, 140, 202, \\78, 119, 109, 248, 38, 15, 193, 238, 31, 83, 252, 131, 
    138, 253, 127, 85, 199, 94, 230, 93, 79, 37, \\87, 40, 111, 8, 161, 249, 17, 
    103, 220, 228, 186, 208, 124, 77, 43, 117, 54, 19, 74, 130, 14, 152, \\102 \}$,
   
\vspace{.02in}
 
        $\{ 266, 42, 81, 178, 100, 257, 260, 258, 273, 15, 243, 160, 21, 119, 76, 
    230, 56, 7, 99, 162, \\164, 170, 25, 57, 153, 156, 54, 9, 19, 182, 92, 39, 
    161, 46, 241, 128, 79, 229, 216, 213, 223, 94,\\ 190, 224, 110, 234, 68, 114, 
    248, 1, 195, 237, 253, 125, 191, 112, 40, 154, 52, 85, 18, 198, \\204, 199, 77
    \}$.   

\vspace{.1in}
    
{\bf LMRH:}

\vspace{.03in}

    $\{ 273, 4, 16, 230, 10, 237, 5, 227, 238, 23, 28, 113, 25, 114, 27, 120, 
    127, 97, 136, 130, 111, \\98, 143, 129, 104, 167, 188, 185, 169, 172, 187, 
    183, 171, 35, 195, 220, 46, 38, 89, 206, 198, 87,\\ 219, 215, 45, 91, 205, 92,
    217, 49, 152, 246, 63, 50, 77, 146, 159, 70, 253, 243, 56, 78, 145, 67, \\254 
    \}$,
   
\vspace{.02in}
 
        $\{ 273, 7, 23, 28, 9, 12, 25, 27, 11, 65, 79, 256, 66, 244, 72, 250, 245, 
    55, 150, 57, 60, 147, 157, \\59, 158, 115, 230, 238, 126, 118, 237, 227, 125, 
    36, 213, 92, 48, 42, 199, 218, 224, 204, 91, 87, \\37, 201, 212, 203, 89, 100,
    165, 182, 112, 106, 134, 170, 176, 131, 189, 179, 101, 141, 164, 142, \\190 \}$.
 
\vspace{.1in}
    
{\bf JOHN:}

\vspace{.03in}

    $\{ 63, 65, 252, 94, 124, 103, 211, 198, 237, 257, 260, 272, 265, 271, 263, 
    262, 267, 129, 178, \\134, 144, 181, 191, 139, 188, 15, 187, 64, 46, 143, 151,
    29, 163, 7, 214, 10, 33, 28, 244, 253, \\195, 219, 48, 21, 62, 229, 236, 206, 
    51, 67, 159, 81, 111, 91, 166, 190, 141, 135, 125, 73, 119,\\ 184, 176, 149, 
    101 \}$,
   
\vspace{.02in}
 
    $\{ 53, 228, 234, 231, 237, 261, 263, 266, 265, 116, 160, 119, 125, 155, 150,
    122, 145, 3, 250, \\8, 61, 52, 253, 195, 206, 58, 55, 200, 14, 201, 9, 247, 
    244, 99, 194, 104, 118, 123, 197, 238, 227, \\113, 128, 233, 110, 232, 105, 
    207, 204, 33, 107, 38, 57, 56, 112, 77, 68, 62, 51, 74, 48, 71, 43, \\102, 97 
    \}$,
   
\vspace{.02in}
 
    $\{ 60, 75, 83, 127, 103, 77, 112, 95, 126, 142, 246, 183, 145, 172, 204, 
    227, 221, 146, 243, \\179, 190, 185, 225, 230, 240, 159, 156, 235, 149, 248, 
    254, 249, 184, 5, 47, 26, 59, 12, 8, 55, \\37, 40, 23, 22, 54, 32, 14, 61, 41,
    133, 195, 140, 192, 168, 201, 229, 209, 177, 169, 253, 157,\\ 239, 219, 247, 
    148 \}$,
   
\vspace{.02in}
 
    $\{ 53, 144, 170, 181, 147, 75, 94, 117, 100, 3, 250, 34, 52, 226, 206, 17, 
    214, 146, 254, 156, \\185, 190, 243, 230, 240, 179, 184, 235, 149, 225, 159, 
    248, 249, 72, 175, 98, 86, 73, 180, 160, \\165, 116, 111, 131, 91, 186, 125, 
    150, 137, 66, 218, 87, 104, 90, 241, 238, 251, 113, 79, 197,\\ 128, 212, 105, 
    207, 232 \}$,
   
\vspace{.02in}
 
$\{ 266, 267, 257, 268, 50, 97, 247, 112, 107, 154, 250, 253, 102, 244, 151, 148,
157, 98, 256, \\111, 108, 146, 241, 246, 101, 251, 159, 156, 149, 114, 203, 127, 
124, 181, 198, 193, 117, 208, \\188, 191, 178, 65, 234, 80, 75, 167, 231, 228, 70,
237, 170, 173, 164, 116, 207, 125, 122, 179, \\194, 197, 119, 204, 190, 185, 184 
\}$,
   
\vspace{.02in}
 
       $\{ 63, 99, 104, 123, 105, 128, 110, 118, 113, 75, 242, 88, 114, 109, 217, 
    208, 231, 194, 207, \\238, 204, 227, 197, 232, 233, 5, 177, 54, 40, 133, 157, 
    23, 169, 131, 230, 189, 145, 180, 199, \\254, 210, 223, 175, 142, 162, 243, 
    202, 235, 160, 67, 159, 81, 111, 91, 166, 190, 141, 135, 125, \\73, 119, 184, 
    176, 149, 101 \}$,
   
\vspace{.02in}
 
    $\{ 63, 66, 76, 90, 69, 84, 79, 93, 87, 212, 215, 256, 221, 218, 251, 246, 
    241, 75, 242, 88, 114, \\109, 217, 208, 231, 5, 177, 54, 40, 133, 157, 23, 
    169, 4, 209, 13, 38, 31, 247, 250, 200, 224, 43,\\ 18, 57, 226, 239, 201, 56, 
    8, 102, 14, 17, 58, 115, 121, 95, 112, 27, 52, 47, 74, 68, 85, 37 \}$,
   
\vspace{.02in}
 
        $\{ 60, 258, 264, 265, 267, 18, 31, 21, 28, 3, 20, 8, 9, 23, 29, 14, 26, 81, 
    105, 86, 98, 111, 110, \\93, 84, 101, 108, 90, 96, 87, 91, 104, 99, 132, 210, 
    135, 166, 171, 213, 196, 205, 161, 176, 199,\\ 141, 202, 138, 223, 220, 131, 
    256, 136, 172, 165, 251, 232, 233, 175, 162, 227, 142, 238, 137, \\241, 246 \}$.

\vspace{.1in}
    
{\bf JOWK:}

\vspace{.03in}

    $\{ 257, 35, 115, 227, 83, 70, 102, 246, 182, 15, 31, 168, 207, 136, 223, 
    152, 56, 1, 30, 170,\\ 205, 17, 137, 210, 14, 148, 155, 138, 221, 59, 52, 169,
    194, 68, 74, 98, 249, 244, 97, 251, 250, \\109, 190, 178, 73, 110, 189, 177, 
    75, 5, 23, 176, 208, 21, 135, 220, 7, 156, 149, 144, 224, 53,\\ 60, 167, 204 
    \}$,
  
\vspace{.02in}
 
        $\{ 273, 4, 203, 202, 10, 9, 201, 11, 196, 17, 209, 18, 30, 29, 210, 222, 
    221, 131, 143, 149, 134, \\156, 136, 160, 151, 52, 59, 162, 58, 57, 174, 161, 
    173, 68, 187, 111, 74, 73, 248, 185, 186, 243, \\75, 104, 99, 102, 180, 255, 
    246, 33, 81, 118, 34, 46, 227, 82, 94, 232, 45, 115, 120, 127, 93, 230, \\239 
    \}$,
  
\vspace{.02in}
 
       $\{ 269, 11, 20, 217, 202, 47, 115, 232, 86, 64, 133, 254, 167, 66, 156, 177,
    109, 1, 17, 57, 205, \\206, 59, 221, 222, 138, 132, 164, 2, 170, 18, 155, 153,
    51, 143, 94, 131, 175, 226, 63, 147, 238,\\ 82, 225, 159, 93, 163, 81, 237, 6,
    24, 41, 198, 200, 123, 216, 214, 122, 36, 116, 8, 42, 22, 43, \\121 \}$,
  
\vspace{.02in}
 
    $\{ 260, 13, 27, 57, 158, 81, 232, 240, 186, 198, 256, 180, 214, 248, 135, 
    167, 82, 5, 28, 109,\\ 63, 123, 147, 124, 101, 41, 47, 78, 67, 3, 32, 253, 60,
    107, 152, 85, 71, 233, 118, 46, 191, 8, 22, \\93, 64, 251, 151, 231, 192, 185,
    88, 238, 246, 12, 21, 77, 51, 43, 159, 111, 115, 121, 76, 110, \\37 \}$,
  
\vspace{.02in}
 
    $\{ 109, 254, 177, 66, 269, 4, 53, 218, 201, 230, 172, 151, 144, 27, 127, 35,
    88, 33, 42, 81, 225, \\186, 84, 235, 121, 113, 180, 187, 185, 34, 122, 82, 
    226, 75, 228, 91, 41, 114, 73, 74, 68, 7, 245, \\224, 197, 30, 252, 247, 256, 
    28, 206, 222, 14, 5, 101, 220, 199, 13, 108, 103, 112, 32, 221, \\205, 29 \}$.

\vspace{.1in}
    
{\bf MATH:}

\vspace{.03in}

$\{ 268, 6, 225, 203, 69, 48, 162, 140, 111, 31, 47, 96, 252, 204, 112, 187, 
    139, 17, 198, \\181, 82, 246, 133, 33, 98, 9, 173, 35, 74, 238, 135, 100, 200,
    3, 196, 127, 68, 232, 28, 131, \\39, 255, 64, 156, 167, 91, 104, 192, 219, 1, 
    86, 24, 66, 230, 4, 21, 177, 231, 87, 243, 165, \\67, 242, 168, 180 \}$,
  
\vspace{.02in}
 
    $\{ 268, 30, 169, 157, 73, 249, 191, 68, 78, 132, 32, 122, 167, 251, 174, 92,
    103, 5, 199, \\87, 246, 226, 253, 211, 36, 102, 43, 180, 56, 208, 17, 26, 129,
    9, 134, 166, 48, 238, 41, 163,\\ 97, 131, 16, 65, 72, 235, 203, 206, 104, 6, 
    256, 10, 100, 225, 88, 195, 27, 189, 186, 237, 40, \\93, 135, 179, 90 \}$,
  
\vspace{.02in}
 
    $\{ 268, 8, 207, 228, 219, 93, 67, 247, 150, 60, 128, 185, 174, 37, 130, 26, 
    97, 12, 191, 51,\\ 27, 81, 121, 245, 230, 41, 58, 106, 176, 99, 66, 36, 116, 
    4, 63, 107, 123, 89, 38, 250, 102, \\164, 83, 50, 163, 10, 34, 127, 249, 16, 
    79, 32, 187, 85, 225, 252, 22, 80, 177, 69, 161, 236,\\ 226, 188, 21 \}$,
  
\vspace{.02in}
 
    $\{ 268, 45, 207, 234, 231, 55, 135, 161, 213, 180, 121, 244, 99, 253, 157, 
    170, 187, 5, 120, \\243, 113, 31, 78, 190, 110, 221, 222, 233, 200, 107, 88, 
    199, 168, 2, 15, 171, 186, 28, 65, 248,\\ 21, 249, 84, 177, 74, 164, 91, 227, 
    238, 10, 203, 111, 85, 20, 73, 215, 209, 22, 115, 117, 105, \\79, 83, 16, 205 
    \}$,
  
\vspace{.02in}
 
    $\{ 268, 22, 114, 46, 41, 241, 195, 155, 149, 105, 206, 175, 201, 128, 76, 
    40, 142, 5, 66, 28, \\136, 226, 168, 15, 165, 77, 99, 196, 255, 236, 39, 67, 
    170, 9, 244, 135, 94, 238, 98, 107, 23, \\164, 185, 120, 100, 144, 147, 133, 
    71, 3, 48, 210, 127, 232, 32, 160, 203, 12, 156, 247, 53, \\123, 20, 251, 239 
    \}$,
  
\vspace{.02in}
 
    $\{ 268, 28, 255, 166, 65, 39, 153, 196, 126, 31, 121, 70, 252, 36, 199, 158,
    161, 20, 127, \\145, 210, 247, 101, 189, 156, 167, 184, 118, 53, 83, 68, 130, 
    90, 15, 32, 220, 57, 236, 84, 42,\\ 251, 253, 146, 63, 222, 117, 26, 183, 205,
    34, 185, 45, 49, 197, 182, 170, 94, 217, 165, 202,\\ 214, 66, 62, 81, 77 \}$,
  
\vspace{.02in}
 
        $\{ 268, 3, 32, 107, 101, 232, 234, 120, 18, 251, 147, 13, 144, 245, 130, 
    157, 122, 27, 174, \\115, 80, 256, 195, 198, 118, 73, 33, 40, 152, 145, 171, 
    249, 30, 2, 163, 106, 196, 229, 79, 203,\\ 9, 72, 48, 172, 141, 238, 39, 134, 
    97, 5, 89, 109, 158, 226, 17, 49, 170, 190, 214, 246, 138, 77, \\121, 37, 194 
    \}$,
  
\vspace{.02in}
 
    $\{ 268, 5, 106, 58, 136, 226, 187, 16, 178, 141, 96, 221, 99, 216, 51, 235, 85, 24,          
    223, 168, \\26, 243, 170, 54, 140, 60, 77, 67, 253, 97, 134, 209, 111,
    1, 6, 242, 80, 230, 191, 176, 18, 225,\\ 92, 21, 171, 188, 95, 75, 245, 10, 
    169, 118, 204, 237, 184, 144, 50, 78, 83, 145, 47, 23, 244, \\107, 213 \}$,
   
\vspace{.02in}
 
   $\{ 268, 21, 142, 253, 86, 242, 37, 205, 105, 177, 190, 26, 42, 102, 194, 89,
    129, 2, 235, 88, \\65, 229, 254, 172, 16, 166, 23, 179, 79, 189, 25, 244, 90, 
    11, 96, 53, 76, 240, 33, 31, 187, 175, \\118, 210, 252, 98, 198, 145, 133, 9, 
    111, 45, 74, 238, 12, 48, 140, 173, 110, 202, 203, 75, 239, \\137, 176 \}$,
  
\vspace{.02in}
 
    $\{ 268, 46, 158, 47, 109, 201, 59, 221, 121, 112, 204, 138, 124, 224, 58, 
    139, 159, 6, 35, 111, \\69, 225, 238, 100, 200, 48, 140, 162, 173, 9, 135, 
    203, 74, 1, 125, 26, 66, 230, 194, 62, 154, 89, \\253, 165, 129, 37, 217, 190,
    102, 2, 235, 113, 65, 229, 64, 172, 16, 50, 150, 166, 127, 219, 79, \\213, 156
    \}$,
  
\vspace{.02in}
 
        $\{ 268, 24, 253, 67, 243, 77, 26, 170, 168, 7, 9, 84, 238, 228, 185, 94, 
    183, 35, 122, 55, 137, \\200, 110, 212, 157, 1, 171, 21, 230, 80, 92, 242, 
    191, 17, 241, 148, 46, 246, 41, 151, 22, 206, 76, \\119, 116, 201, 172, 79, 
    175, 8, 13, 198, 124, 227, 113, 207, 234, 150, 93, 33, 44, 159, 88, 179, \\186
    \}$,
  
\vspace{.02in}
 
    $\{ 268, 43, 129, 102, 220, 208, 149, 114, 63, 9, 72, 250, 238, 84, 163, 183,
    29, 40, 42, 152, \\195, 205, 125, 115, 154, 55, 110, 135, 57, 212, 222, 100, 
    137, 10, 144, 112, 210, 237, 88, 234, \\107, 179, 83, 139, 13, 53, 213, 50, 
    184, 5, 41, 36, 158, 226, 85, 235, 206, 178, 96, 199, 16, 121, \\116, 151, 187
    \}$,
  
\vspace{.02in}
 
    $\{ 268, 4, 9, 84, 89, 231, 238, 183, 190, 32, 162, 154, 69, 251, 40, 125, 
    195, 54, 102, 129, 111, \\209, 140, 63, 220, 1, 100, 92, 230, 57, 135, 191, 
    222, 33, 193, 107, 124, 198, 112, 156, 38, 53, \\50, 144, 159, 139, 127, 210, 
    213, 5, 11, 196, 96, 226, 116, 82, 240, 41, 153, 39, 187, 151, \\181, 206, 126
    \}$,
  
\vspace{.02in}
 
    $\{ 268, 12, 52, 81, 239, 105, 215, 182, 142, 49, 97, 108, 214, 60, 134, 143,
    223, 21, 171, 202, \\80, 242, 120, 45, 147, 7, 14, 180, 94, 228, 87, 233, 185,
    33, 193, 107, 124, 198, 112, 156, 38, 53,\\ 50, 144, 159, 139, 127, 210, 213, 
    5, 11, 196, 96, 226, 116, 82, 240, 41, 153, 39, 187, 151, 181, \\206, 126 \}$,
  
\vspace{.02in}
 
    $\{ 268, 35, 135, 100, 202, 200, 110, 137, 45, 11, 119, 76, 240, 56, 148, 
    175, 211, 20, 139, 83, \\204, 247, 184, 112, 47, 19, 248, 84, 94, 183, 185, 
    29, 250, 3, 41, 68, 101, 232, 106, 235, 206, 38, \\167, 16, 130, 172, 141, 79,
    193, 14, 191, 77, 26, 233, 256, 236, 92, 89, 170, 15, 253, 171, 27, 80, \\190 
    \}$,
  
\vspace{.02in}
 
    $\{ 268, 30, 98, 93, 33, 249, 186, 133, 198, 5, 158, 70, 226, 221, 121, 161, 
    58, 18, 81, 251, 245, \\188, 182, 95, 32, 34, 97, 111, 197, 48, 134, 140, 203,
    3, 41, 68, 101, 232, 106, 235, 206, 38, 167, \\16, 130, 172, 141, 79, 193, 14,
    191, 77, 26, 233, 256, 236, 92, 89, 170, 15, 253, 171, 27, 80, 190 \}$.

\vspace{.1in}
    
{\bf SEMI2:}

\vspace{.03in}

  $\{ 273, 81, 97, 139, 98, 85, 206, 96, 101, 137, 199, 131, 82, 205, 112, 204, 134, 3, 14, 73, 147, \\11, 
    71, 158, 12, 233, 231, 70, 156, 237, 155, 77, 230, 17, 222, 245, 51, 21, 192, 176, 220, 35, 124,\\ 241, 
    162, 126, 59, 178, 43, 22, 209, 249, 50, 25, 183, 167, 213, 47, 116, 246, 173, 120, 64, 189, \\42 \}$,
  
\vspace{.02in}
 
        $\{ 273, 87, 135, 197, 88, 94, 105, 136, 142, 98, 195, 194, 90, 99, 138, 101, 201, 2, 74, 230, 9, \\5,
    148, 78, 72, 157, 240, 235, 3, 159, 71, 156, 225, 34, 173, 220, 41, 37, 185, 164, 172, 178, 223, \\221, 
    35, 179, 175, 181, 212, 20, 254, 125, 29, 31, 55, 250, 247, 56, 116, 124, 28, 62, 248, 58, 127 \}$,
  
\vspace{.02in}
 
    $\{ 273, 38, 75, 110, 48, 47, 115, 65, 68, 120, 106, 98, 44, 117, 77, 119, 105, 86, 245, 184, 96, \\95, 
    230, 243, 247, 239, 183, 179, 92, 240, 248, 236, 181, 19, 139, 152, 21, 24, 209, 129, 132, 221, \\151, 
    147, 23, 219, 141, 212, 149, 6, 176, 61, 16, 15, 193, 166, 172, 205, 52, 49, 12, 203, 175, 196, \\59 \}$,
  
\vspace{.02in}
 
        $\{ 273, 35, 251, 152, 196, 45, 198, 156, 247, 255, 145, 242, 197, 48, 206, 42, 153, 1, 6, 190, \\116, 8,
    184, 124, 14, 136, 182, 129, 117, 142, 177, 134, 121, 20, 210, 101, 226, 21, 76, 91, 223, \\171, 100, 
    167, 239, 57, 73, 60, 87, 23, 220, 107, 230, 27, 66, 88, 217, 165, 103, 164, 238, 63, 79, \\50, 81 \}$,
  
\vspace{.02in}
 
    $\{ 273, 26, 213, 255, 127, 32, 178, 48, 212, 42, 51, 242, 191, 61, 114, 172, 169, 25, 211, 246,\\ 118, 
    28, 185, 38, 221, 46, 56, 254, 188, 49, 126, 175, 162, 4, 12, 69, 226, 5, 73, 227, 9, 237, 155, \\68, 
    76, 151, 239, 160, 154, 84, 109, 131, 196, 85, 206, 138, 99, 144, 106, 141, 198, 112, 197, 86,\\ 94 \}$,
  
\vspace{.02in}
 
    $\{ 273, 49, 175, 153, 135, 56, 130, 145, 162, 152, 171, 156, 143, 167, 139, 57, 60, 38, 246, 81,\\ 236, 
    46, 80, 103, 254, 107, 118, 88, 74, 126, 233, 185, 188, 34, 250, 83, 234, 47, 78, 97, 256, 104,\\ 114, 
    93, 70, 127, 240, 191, 178, 1, 6, 32, 217, 8, 27, 213, 14, 212, 202, 26, 23, 208, 220, 199, \\203 \}$,
   
\vspace{.02in}
 
       $\{ 273, 2, 16, 132, 10, 136, 15, 133, 129, 146, 197, 160, 154, 200, 193, 159, 196, 113, 185, 117, \\116,
    182, 190, 120, 188, 33, 254, 246, 37, 36, 252, 249, 40, 66, 232, 94, 80, 74, 109, 228, 229, 107,\\ 92, 
    89, 79, 103, 225, 99, 86, 22, 224, 174, 25, 28, 63, 210, 223, 64, 172, 169, 30, 58, 218, 50, \\166 \}$,
  
\vspace{.02in}
 
    $\{ 273, 1, 131, 5, 4, 139, 135, 8, 141, 118, 180, 181, 121, 124, 184, 177, 126, 33, 254, 252, 37, \\36, 
    249, 246, 40, 147, 157, 196, 155, 151, 197, 193, 200, 19, 209, 169, 27, 23, 58, 213, 212, 63, \\166, 
    174, 29, 50, 216, 64, 172, 67, 230, 87, 75, 71, 102, 233, 236, 105, 93, 83, 77, 108, 238, 110,\\ 91 \}$,
  
\vspace{.02in}
 
    $\{ 273, 114, 179, 187, 128, 122, 183, 127, 189, 1, 136, 5, 4, 129, 8, 132, 133, 33, 254, 37, 36, \\40, 
    252, 249, 246, 150, 156, 204, 153, 198, 158, 206, 201, 17, 215, 173, 21, 20, 56, 221, 211, 49,\\ 24, 
    167, 171, 163, 219, 52, 53, 67, 230, 86, 75, 71, 107, 233, 236, 103, 77, 89, 92, 94, 238, 99, \\109 \}$,
  
\vspace{.02in}
 
    $\{ 273, 35, 163, 218, 46, 44, 191, 174, 172, 177, 213, 209, 43, 181, 171, 186, 223, 81, 131, 199, \\95, 
    90, 104, 142, 140, 98, 201, 198, 85, 112, 139, 100, 205, 18, 241, 122, 24, 20, 56, 255, 250, 50,\\ 117, 
    113, 32, 64, 245, 52, 127, 2, 67, 228, 8, 4, 158, 78, 76, 147, 240, 226, 16, 155, 75, 156, 232 \}$,
  
\vspace{.02in}
 
    $\{ 273, 35, 72, 109, 37, 40, 116, 71, 67, 123, 100, 97, 39, 125, 69, 113, 107, 3, 174, 62, 5, 8, 199, \\
    170, 162, 197, 58, 50, 7, 200, 169, 195, 57, 83, 253, 191, 85, 88, 234, 244, 241, 233, 188, 182, 87, \\
    238, 251, 226, 192, 17, 142, 157, 27, 29, 218, 138, 130, 217, 148, 145, 20, 222, 137, 210, 155 \}$,
  
\vspace{.02in}
 
    $\{ 273, 20, 246, 125, 29, 31, 59, 256, 251, 49, 116, 124, 28, 54, 241, 64, 127, 36, 172, 211, 45,\\ 47, 
    179, 175, 173, 181, 213, 217, 44, 185, 164, 178, 210, 84, 133, 203, 93, 95, 106, 131, 130, 110, \\193, 
    198, 92, 104, 137, 103, 208, 1, 69, 227, 11, 16, 160, 67, 66, 150, 229, 233, 6, 145, 73, 155, \\226 \}$,
  
\vspace{.02in}
 
    $\{ 273, 6, 181, 78, 14, 11, 210, 180, 178, 70, 71, 7, 223, 213, 191, 75, 212, 51, 92, 118, 61, 57, \\
    134, 89, 93, 126, 123, 60, 142, 139, 83, 119, 135, 17, 101, 249, 24, 32, 199, 100, 98, 252, 243, 26, \\
    203, 206, 111, 253, 198, 35, 149, 162, 45, 41, 237, 148, 146, 175, 165, 44, 227, 236, 159, 164, \\233 
    \}$,
  
\vspace{.02in}
 
    $\{ 273, 22, 98, 245, 30, 27, 206, 111, 101, 244, 242, 23, 198, 199, 100, 255, 203, 6, 179, 78, 14, \\11,
    217, 189, 185, 70, 71, 7, 220, 211, 188, 75, 221, 35, 147, 162, 45, 41, 228, 157, 153, 175, 165, \\44, 
    229, 239, 156, 164, 226, 50, 93, 113, 63, 53, 139, 83, 92, 120, 128, 52, 135, 134, 89, 122, 142 \}$,
  
\vspace{.02in}
 
    $\{ 273, 35, 147, 166, 45, 41, 238, 157, 153, 44, 174, 171, 156, 230, 231, 167, 235, 3, 184, 67, \\13, 9,
    209, 177, 186, 12, 77, 73, 192, 216, 224, 76, 218, 50, 88, 126, 63, 53, 141, 81, 90, 52, 118, \\119, 96,
    131, 140, 123, 137, 22, 112, 241, 30, 27, 199, 106, 97, 23, 248, 256, 104, 203, 206, 250, \\198 \}$,
  
\vspace{.02in}
 
    $\{ 273, 34, 145, 161, 47, 37, 239, 152, 160, 36, 168, 176, 154, 226, 228, 170, 229, 54, 94, 127, \\62, 
    59, 143, 86, 87, 55, 114, 116, 91, 130, 132, 117, 133, 1, 190, 65, 8, 16, 214, 182, 183, 10, 72, \\80, 
    187, 222, 219, 74, 215, 19, 112, 242, 29, 25, 199, 106, 97, 28, 255, 245, 104, 203, 206, 244, \\198 \}$,
  
\vspace{.02in}
 
    $\{ 273, 33, 156, 172, 40, 48, 42, 240, 153, 157, 169, 173, 147, 163, 234, 225, 232, 54, 87, 128,\\ 62, 
    59, 55, 138, 91, 94, 122, 113, 86, 120, 144, 136, 129, 3, 187, 73, 13, 9, 12, 219, 183, 182, 76, \\67, 
    190, 77, 215, 214, 222, 17, 110, 255, 24, 32, 26, 201, 102, 103, 242, 244, 107, 245, 204, 195,\\ 205 \}$,
  
\vspace{.02in}
 
    $\{ 273, 34, 159, 167, 47, 37, 36, 230, 146, 148, 171, 174, 149, 166, 238, 235, 231, 6, 177, 75, \\14, 
    11, 7, 209, 184, 192, 71, 70, 186, 78, 216, 224, 218, 18, 106, 248, 31, 21, 20, 194, 112, 104, \\241, 
    250, 97, 256, 207, 197, 196, 54, 95, 128, 62, 59, 55, 141, 82, 84, 122, 113, 85, 120, 131, \\140, 137 \}$,
  
\vspace{.02in}
 
    $\{ 273, 17, 107, 235, 25, 28, 168, 98, 111, 172, 226, 239, 24, 169, 103, 161, 231, 2, 192, 117, \\11, 7,
    136, 179, 189, 140, 118, 126, 15, 137, 186, 129, 116, 36, 158, 200, 46, 38, 249, 148, 149, \\241, 204, 
    201, 37, 248, 150, 252, 193, 52, 87, 70, 62, 54, 223, 95, 82, 215, 69, 68, 53, 219, 91, \\210, 78 \}$,
  
\vspace{.02in}
 
    $\{ 273, 82, 136, 194, 91, 87, 108, 140, 137, 97, 203, 199, 95, 104, 129, 105, 207, 20, 242, 115, \\30, 
    22, 56, 251, 247, 57, 128, 122, 21, 60, 255, 49, 125, 3, 75, 239, 16, 10, 157, 66, 79, 160, 231, \\235, 
    13, 154, 71, 147, 226, 35, 163, 217, 48, 42, 177, 176, 170, 188, 209, 216, 45, 185, 173, 184, \\220 \}$,
  
\vspace{.02in}
 
    $\{ 273, 17, 187, 147, 25, 28, 84, 178, 191, 94, 160, 154, 24, 86, 183, 85, 157, 3, 111, 59, 16, \\10, 
    48, 103, 107, 35, 50, 63, 13, 45, 98, 42, 55, 65, 201, 119, 73, 76, 162, 193, 200, 171, 127, \\114, 72, 
    167, 204, 175, 123, 131, 225, 224, 144, 138, 252, 233, 236, 248, 211, 221, 141, 241, \\232, 249, 218 \}$.

\vspace{.1in}
    
{\bf SEMI4:}

\vspace{.03in}

   $\{ 273, 49, 55, 80, 163, 56, 78, 164, 59, 70, 112, 74, 110, 106, 173, 102, 165,
    9, 16, 119, 197, \\12, 116, 203, 10, 117, 44, 123, 42, 41, 196, 48, 199, 17, 149,
    232, 215, 24, 241, 142, 148, 248, \\85, 225, 186, 84, 219, 192, 134, 25, 153, 
    236, 223, 28, 246, 136, 156, 254, 86, 233, 183, 94, 210, \\187, 129 \}$,
  
\vspace{.02in}
 
        $\{ 273, 17, 106, 49, 20, 29, 27, 82, 98, 105, 52, 61, 110, 59, 90, 94, 89, 70, 
    209, 237, 76, 79,\\ 80, 176, 212, 221, 235, 225, 219, 228, 175, 172, 166, 6, 182,
    158, 12, 15, 16, 41, 188, 191, 153, \\146, 192, 154, 46, 42, 34, 115, 129, 255, 
    119, 120, 117, 201, 132, 141, 256, 246, 139, 252, 206,\\ 202, 194 \}$,
  
\vspace{.02in}
 
    $\{ 273, 164, 238, 173, 166, 234, 229, 176, 227, 68, 210, 77, 70, 217, 215, 80, 
    216, 36, 153, 152, \\45, 38, 48, 146, 151, 3, 10, 180, 5, 182, 14, 192, 189, 18, 
    113, 99, 25, 23, 132, 123, 124, 141, 101, \\106, 24, 134, 127, 144, 110, 50, 195,
    87, 57, 55, 245, 197, 202, 243, 88, 82, 56, 254, 206, 250, \\89 \}$,
  
\vspace{.02in}
 
    $\{ 273, 3, 7, 115, 13, 119, 11, 125, 123, 35, 208, 202, 39, 45, 207, 43, 194, 
    50, 79, 58, 63, 74, \\64, 80, 66, 102, 161, 165, 108, 110, 164, 168, 105, 83, 
    137, 217, 87, 93, 183, 142, 140, 189, 91, \\222, 220, 214, 134, 179, 187, 22, 
    146, 239, 28, 30, 246, 154, 159, 249, 25, 240, 226, 234, 160, \\252, 254 \}$,
   
\vspace{.02in}
 
   $\{ 273, 193, 221, 224, 206, 204, 197, 210, 215, 2, 13, 32, 7, 23, 16, 29, 18, 
    132, 170, 171, 137,\\ 134, 163, 175, 136, 49, 53, 145, 62, 60, 149, 156, 158, 34,
    184, 88, 45, 39, 103, 182, 185, 112, 86, \\89, 48, 98, 180, 109, 84, 68, 117, 
    241, 73, 70, 228, 124, 126, 233, 254, 252, 72, 230, 113, 232, \\245 \}$,
  
\vspace{.02in}
 
        $\{ 273, 133, 186, 159, 136, 143, 242, 178, 177, 244, 160, 149, 144, 250, 180, 
    241, 152, 5, 120, \\195, 8, 15, 44, 117, 128, 35, 199, 204, 16, 38, 127, 39, 198,
    51, 79, 169, 55, 60, 106, 80, 69, 97,\\ 174, 173, 54, 98, 72, 100, 171, 19, 225, 
    210, 23, 28, 83, 228, 234, 92, 218, 212, 22, 87, 226, 86,\\ 209 \}$,
  
\vspace{.02in}
 
        $\{ 273, 259, 270, 266, 261, 146, 256, 157, 149, 240, 254, 247, 156, 225, 
    238, 241, 231, 116, \\174, 123, 117, 190, 175, 162, 119, 182, 191, 166, 178, 
    50, 105, 64, 51, 82, 97, 102, 59, 83, 91, \\101, 96, 37, 140, 40, 48, 71, 131,
    134, 47, 76, 70, 135, 67, 3, 207, 10, 5, 223, 193, 203, 14, \\220, 209, 204, 
    219 \}$,
  
\vspace{.02in}
 
    $\{ 270, 21, 118, 186, 173, 8, 256, 236, 98, 89, 51, 212, 139, 193, 151, 46, 
    79, 9, 247, 225, 103,\\ 86, 63, 217, 143, 197, 157, 33, 74, 14, 246, 228, 99, 
    83, 53, 221, 144, 208, 158, 42, 68, 3, 250, \\240, 110, 96, 61, 213, 132, 195,
    148, 38, 78, 15, 255, 231, 97, 93, 55, 209, 137, 202, 150, 41, \\69 \}$.

\begin{Note}
Magma \cite{magma} was used for all of the computations performed in this paper.
\end{Note}

\end{document}